\documentclass[12pt]{article}
\usepackage{amsfonts}
\usepackage{full page,
amssymb, amscd, amsmath,graphicx, 
}

\title{{\bf On associative algebras, modules and twisted modules for
    vertex operator algebras}} \author{Jinwei Yang} \date{}
    \begin{document}
    \bibliographystyle{alpha}

\newtheorem{thm}{Theorem}[section]
\newtheorem{defn}[thm]{Definition}
\newtheorem{prop}[thm]{Proposition}
\newtheorem{cor}[thm]{Corollary}
\newtheorem{lemma}[thm]{Lemma}
\newtheorem{rema}[thm]{Remark}
\newtheorem{app}[thm]{Application}
\newtheorem{prob}[thm]{Problem}
\newtheorem{conv}[thm]{Convention}
\newtheorem{conj}[thm]{Conjecture}
\newtheorem{cond}[thm]{Condition}
    \newtheorem{exam}[thm]{Example}
\newtheorem{assum}[thm]{Assumption}
     \newtheorem{nota}[thm]{Notation}
\newcommand{\halmos}{\rule{1ex}{1.4ex}}
\newcommand{\pfbox}{\hspace*{\fill}\mbox{$\halmos$}}

\newcommand{\nn}{\nonumber \\}

 \newcommand{\res}{\mbox{\rm Res}}
 \newcommand{\ord}{\mbox{\rm ord}}
\renewcommand{\hom}{\mbox{\rm Hom}}
\newcommand{\edo}{\mbox{\rm End}\ }
 \newcommand{\pf}{{\it Proof.}\hspace{2ex}}
 \newcommand{\epf}{\hspace*{\fill}\mbox{$\halmos$}}
 \newcommand{\epfv}{\hspace*{\fill}\mbox{$\halmos$}\vspace{1em}}
 \newcommand{\epfe}{\hspace{2em}\halmos}
\newcommand{\nord}{\mbox{\scriptsize ${\circ\atop\circ}$}}
\newcommand{\wt}{\mbox{\rm wt}\ }
\newcommand{\swt}{\mbox{\rm {\scriptsize wt}}\ }
\newcommand{\lwt}{\mbox{\rm wt}^{L}\;}
\newcommand{\rwt}{\mbox{\rm wt}^{R}\;}
\newcommand{\slwt}{\mbox{\rm {\scriptsize wt}}^{L}\,}
\newcommand{\srwt}{\mbox{\rm {\scriptsize wt}}^{R}\,}
\newcommand{\clr}{\mbox{\rm clr}\ }
\newcommand{\tr}{\mbox{\rm Tr}}
\newcommand{\C}{\mathbb{C}}
\newcommand{\Z}{\mathbb{Z}}
\newcommand{\R}{\mathbb{R}}
\newcommand{\Q}{\mathbb{Q}}
\newcommand{\N}{\mathbb{N}}
\newcommand{\CN}{\mathcal{N}}
\newcommand{\F}{\mathcal{F}}
\newcommand{\I}{\mathcal{I}}
\newcommand{\V}{\mathcal{V}}
\newcommand{\one}{\mathbf{1}}
\newcommand{\BY}{\mathbb{Y}}
\newcommand{\ds}{\displaystyle}

        \newcommand{\ba}{\begin{array}}
        \newcommand{\ea}{\end{array}}
        \newcommand{\be}{\begin{equation}}
        \newcommand{\ee}{\end{equation}}
        \newcommand{\bea}{\begin{eqnarray}}
        \newcommand{\eea}{\end{eqnarray}}
         \newcommand{\lbar}{\bigg\vert}
        \newcommand{\p}{\partial}
        \newcommand{\dps}{\displaystyle}
        \newcommand{\bra}{\langle}
        \newcommand{\ket}{\rangle}

        \newcommand{\ob}{{\rm ob}\,}
        \renewcommand{\hom}{{\rm Hom}}

\newcommand{\A}{\mathcal{A}}
\newcommand{\Y}{\mathcal{Y}}
    \maketitle

\newcommand{\dlt}[3]{#1 ^{-1}\delta \bigg( \frac{#2 #3 }{#1 }\bigg) }

\newcommand{\dlti}[3]{#1 \delta \bigg( \frac{#2 #3 }{#1 ^{-1}}\bigg) }

 \makeatletter
\newlength{\@pxlwd} \newlength{\@rulewd} \newlength{\@pxlht}
\catcode`.=\active \catcode`B=\active \catcode`:=\active
\catcode`|=\active
\def\sprite#1(#2,#3)[#4,#5]{
   \edef\@sprbox{\expandafter\@cdr\string#1\@nil @box}
   \expandafter\newsavebox\csname\@sprbox\endcsname
   \edef#1{\expandafter\usebox\csname\@sprbox\endcsname}
   \expandafter\setbox\csname\@sprbox\endcsname =\hbox\bgroup
   \vbox\bgroup
  \catcode`.=\active\catcode`B=\active\catcode`:=\active\catcode`|=\active
      \@pxlwd=#4 \divide\@pxlwd by #3 \@rulewd=\@pxlwd
      \@pxlht=#5 \divide\@pxlht by #2
      \def .{\hskip \@pxlwd \ignorespaces}
      \def B{\@ifnextchar B{\advance\@rulewd by \@pxlwd}{\vrule
         height \@pxlht width \@rulewd depth 0 pt \@rulewd=\@pxlwd}}
      \def :{\hbox\bgroup\vrule height \@pxlht width 0pt depth
0pt\ignorespaces}
      \def |{\vrule height \@pxlht width 0pt depth 0pt\egroup
         \prevdepth= -1000 pt}
   }
\def\endsprite{\egroup\egroup}
\catcode`.=12 \catcode`B=11 \catcode`:=12 \catcode`|=12\relax
\makeatother

\def\hboxtr{\FormOfHboxtr} 
\sprite{\FormOfHboxtr}(25,25)[0.5 em, 1.2 ex] 

:BBBBBBBBBBBBBBBBBBBBBBBBB | :BB......................B |
:B.B.....................B | :B..B....................B |
:B...B...................B | :B....B..................B |
:B.....B.................B | :B......B................B |
:B.......B...............B | :B........B..............B |
:B.........B.............B | :B..........B............B |
:B...........B...........B | :B............B..........B |
:B.............B.........B | :B..............B........B |
:B...............B.......B | :B................B......B |
:B.................B.....B | :B..................B....B |
:B...................B...B | :B....................B..B |
:B.....................B.B | :B......................BB |
:BBBBBBBBBBBBBBBBBBBBBBBBB |

\endsprite
\def\shboxtr{\FormOfShboxtr} 
\sprite{\FormOfShboxtr}(25,25)[0.3 em, 0.72 ex] 

:BBBBBBBBBBBBBBBBBBBBBBBBB | :BB......................B |
:B.B.....................B | :B..B....................B |
:B...B...................B | :B....B..................B |
:B.....B.................B | :B......B................B |
:B.......B...............B | :B........B..............B |
:B.........B.............B | :B..........B............B |
:B...........B...........B | :B............B..........B |
:B.............B.........B | :B..............B........B |
:B...............B.......B | :B................B......B |
:B.................B.....B | :B..................B....B |
:B...................B...B | :B....................B..B |
:B.....................B.B | :B......................BB |
:BBBBBBBBBBBBBBBBBBBBBBBBB |

\endsprite

\vspace{2em}

\begin{abstract}
  We give a new construction of functors from the category of modules
  for the associative algebras $A_n(V)$ and $A_g(V)$ associated with a
  vertex operator algebra $V$, defined by Dong, Li and Mason, to the
  category of admissible $V$-modules and admissible twisted
  $V$-modules, respectively, using the method developed in the joint work
  \cite{HY1} with Y.-Z. Huang. The functors were first constructed by Dong, Li and
  Mason, but the importance of the new method, as in \cite{HY1}, is
  that we can apply the method to study objects without the commutator
  formula in the representation theory of vertex operator algebras.
\end{abstract}



\renewcommand{\theequation}{\thesection.\arabic{equation}}
\renewcommand{\thethm}{\thesection.\arabic{thm}}
\setcounter{equation}{0} \setcounter{thm}{0} 
\date{}
\maketitle

\setcounter{equation}{0}

\section{Introduction}
This paper is a continuation of the paper \cite{HY1}. The aim is to
prove results in the representation theory of vertex operator
algebras, in particular, for modules, without using the standard commutator
formula. The commutator formula for vertex operators plays a very
important role in the representation theory of vertex operator
algebras since it allows one to apply many techniques in
Lie algebra representation theory to study vertex operator algebras
and their modules. However, for some important objects in the
representation theory of vertex operator algebras such as intertwining
operators (or more generally, logarithmic intertwining operators),
there is no commutator formula for two intertwining operators and therefore we have to use the
associativity of intertwining operators.

In \cite{HY1}, jointly with Y.-Z. Huang, the author gave a formula for the
residues of certain formal series involving iterates of vertex
operators obtained using the weak associativity and the lower
truncation property of vertex operators. We proved that the weak
associativity for an admissible module is equivalent to this residue
formula together with a formula that expresses products of components
of vertex operators as linear combinations of iterates of components
of vertex operators given in \cite{DLM1} and \cite{L}.  We applied
this result to give a new construction of admissible modules for an
$\N$-graded vertex algebra $V$ from modules for its Zhu algebra $A(V)$.

In this paper, we use the method in \cite{HY1}, but in more general
settings, to construct a functor from the category of modules for the
associative algebra $A_n(V)$, defined in \cite{DLM1}, generalizing the Zhu algebra associated with
a vertex operator algebra $V$ for $n \in \N$, to
the category of admissible $V$-modules. We also use the method in
\cite{HY1} to construct a functor from the category of modules for the
``twisted" generalization of the Zhu algebra $A_g(V)$, defined in \cite{DLM2}, associated with a
vertex operator algebra $V$ and a finite order automorphism $g$ of $V$,
to the category of admissible $g$-twisted
$V$-modules.

The associative algebra $A_n(V)$ plays an important
role in the representation theory of vertex operator
algebras. One example is the study of logarithmic intertwining
operators among generalized modules for a vertex operator algebra. In
\cite{HY2}, jointly with Y.-Z. Huang, the author proved that the space of
logarithmic intertwining operators among generalized modules is naturally
isomorphic to the space of homomorphisms between suitable modules for
$A_n(V)$.

The twisted generalization of the Zhu algebra $A_g(V)$ was introduced to
study twisted modules for a vertex operator algebra with an
automorphism of finite order. In \cite{H}, Y.-Z. Huang generalized the
notion of twisted module for a vertex operator algebra with a
finite-order automorphism to the notion of generalized twisted module
for a vertex operator algebra with an automorphism of {\it not
  necessarily finite order}, using logarithmic conformal field
theory. It is natural to define a suitable associative algebra, generalizing the Zhu algebra, associated with a non-finite-order automorphism of $V$, and to construct generalized twisted
modules in the sense of \cite{H} from certain modules for that associative algebra. For
these generalized twisted modules, the twisted vertex operators
involve the logarithm of the variable and thus do not have a
commutator formula. This motivates us in the present paper to
discover a new construction of twisted $V$-modules from
$A_g(V)$-modules without using commutator formula.

The formula for the residues of certain formal series involving iterates of
vertex operators, discovered in \cite{HY1}, has two undetermined parameters satisfying the lower truncation property. By
specializing the parameters to suitable values in the residue formula, we obtain the actions of
the relations used to define associative algebras $A_n(V)$, as well as
$A_g(V)$, on the admissible $V$-module. This coincidence motivates us
to discover the fact that the relations for defining various
generalizations of the Zhu algebra are special cases of the residue formula, and hence
are implied by the weak associativity, straightforwardly. Based on this fact, we construct functors from the category of modules for the associative algebra to the categories of suitable modules and twisted modules for
vertex operator algebras. The most technical part that we prove in
this paper is the converse of the fact, that is, that the relations for
defining various generalizations of the Zhu algebra imply the residue
formula, provided that the lower truncation property and the formula
that expresses products of components of vertex operators as linear
combinations of iterates of components of vertex operators, mentioned
above, hold.

The functors we have constructed in this paper satisfy the same
universal property as the functors constructed in \cite{DLM1} and
\cite{DLM2}. The importance of our construction is that it allows us
to use the method in the present paper to give constructions and prove
results in the cases where there is no commutator formula.

There are various other generalizations of the Zhu algebra associated with
a vertex operator algebra in the literature, such as those in \cite{DLM3},
\cite{MT} and \cite{VE}. The actions of the relations used to define
these associative algebras on certain $V$-modules can be obtained from
the residue formula by specializing the parameters to suitable
values. It is expected that we can use the method in \cite{HY1} and the current paper to
give new constructions of functors between module categories for these
associative algebras and suitable module categories for the vertex
operator algebra.

This paper is organized as follows: We recall the main theorem of
\cite{HY1} and deduce some corollaries to motivate this work in
Section 2. In Section 3, we recall definitions and properties of
the associative algebras $A_n(V)$ and $A_g(V)$. In Section 4 and 5, we
use the main theorem to construct a functor from the category of
$A_n(V)$-modules to the category of admissible $V$-modules and prove a
universal property of this functor. In Section 6, we apply the main
theorem to the twisted module case and construct a functor from the
category of $A_g(V)$-modules to the category of admissible $g$-twisted
$V$-modules.

\paragraph{Acknowledgments}
I would like to thank Prof. Yi-Zhi Huang for inspiring me to think
about this direction and thank Profs. James Lepowsky, Haisheng Li and Katrina Barron for helpful suggestions. I also would like to express my gratitude to Prof. Ping Li for his support during the year 2013--2014.

\setcounter{equation}{0}
\section{An equivalent condition for the associativity}
Throughout this paper, we will let $(V, Y, {\bf 1}, \omega)$ denote a
vertex operator algebra. We state the following theorem from
\cite{HY1}, using a more general setting:
\begin{thm}\label{main theorem}
Let $V$ be a vertex operator algebra, $W$ a vector space and
$Y_W$ a linear map from $V\otimes W$ to $W((x))$. For $v\in V$ and $w\in W$,
as for a $V$-module, we denote the image of $u\otimes w$ under $Y_W$
by $Y_W(u, x)w$ and $\res_{x}x^{n}Y_W(u, x)w$ by $u_{n}w$.
Let $u, v\in V$ and $w\in W$ and let  $k, l \in \R$ such that
\begin{equation}\label{eq6}
v_nw = 0\;\;\; {\rm for}\; n \geq k
\end{equation}
and
\begin{equation}\label{eq7}
u_nw = 0\;\;\;{\rm for}\; n \geq l.
\end{equation}
Then the weak associativity for $W$ in the sense that
\begin{equation}\label{wk-assoc}
(x_{0}+x_{2})^{l}Y_W(u, x_0 + x_2)Y_W(v, x_2)w
=(x_{2}+x_{0})^{l}Y_W(Y(u, x_0)v, x_2)w
\end{equation}
is equivalent to the following two properties: For any $p \in l + \Z$, $q \in k + \Z$,
\begin{eqnarray}\label{eq1-1}
\lefteqn{\res_{x_{0}}\res_{x_{2}}(x_{0}+x_{2})^{p}x_{2}^{q}Y_W(u, x_{0}+x_{2})
Y_W(v, x_{2})w}\nn
&&=\res_{x_{0}}\res_{x_{2}}f(x_{0}, x_{2})x_{2}^{q}(x_{2}+x_{0})^{l}
Y_W(Y(u, x_{0})v, x_{2})w,
\end{eqnarray}
where
\begin{equation}\label{defp}
f(x_{0}, x_{2})=\sum_{i=0}^{k-q-1}\binom{p-l}{i}x_{0}^{p-l-i}x_{2}^{i},
\end{equation}
and
\begin{equation}\label{upper-truncation-3}
\res_{x_{0}}\res_{x_{2}}x_{0}^{p-l-i}x_{2}^{q+i}(x_{2}+x_{0})^{l}Y_{W}(Y(u, x_{0})v, x_{2})
w=0
\end{equation}
for $i\ge k-q$.
\end{thm}

\pf The statement is slightly different from Theorem 3.1 in \cite{HY1}
because here $k, l$ can be real numbers instead of only integers. The
proof for ``only if" part is the same as \cite{HY1}, we give a proof for
the ``if" part here.

The Laurent polynomial $p(x_{0}, x_{2})$
is in fact the first $k-q$ terms of the formal series $(x_{0}+x_{2})^{p-l}$.
But from (\ref{upper-truncation-3}), we obtain
\begin{equation}\label{upper-truncation-7}
\res_{x_{0}}\res_{x_{2}}((x_{0}+x_{2})^{p-l}-f(x_{0}, x_{2}))x_{2}^{q}
((x_{2}+x_{0})^{l}Y_{W}(Y(u, x_{0})v, x_{2}))
w=0.
\end{equation}
Combining (\ref{eq1-1}) and (\ref{upper-truncation-7}),
we obtain
\begin{eqnarray*}
\lefteqn{\res_{x_{0}}\res_{x_{2}}(x_{0}+x_{2})^{p-l}x_{2}^{q}((x_{0}+x_{2})^{l}
Y_{W}(u, x_{0}+x_{2})
Y_{W}(v, x_{2}))w}\nn
&&=\res_{x_{0}}\res_{x_{2}}(x_{0}+x_{2})^{p-l}x_{2}^{q}((x_{2}+x_{0})^{l}
Y_{W}(Y(u, x_{0})v, x_{2}))w.
\end{eqnarray*}

On the other hand, we have
\begin{equation}\label{upper-truncation-5}
\res_{x_{0}}\res_{x_{2}}x_{0}^{p-l-i}x_{2}^{q+i}(x_{0}+x_{2})^{l}Y_{W}(u, x_{0}+x_{2})
Y_{W}(v, x_{2})
w=0
\end{equation}
for $i\ge k-q$. From (\ref{upper-truncation-3}) and (\ref{upper-truncation-5}),
we obtain
\begin{eqnarray}\label{upper-truncation-6}
\lefteqn{\res_{x_{0}}\res_{x_{2}}x_{0}^{p-l-i}x_{2}^{q+i}(x_{0}+x_{2})^{l}Y_{W}(u, x_{0}+x_{2})
    Y_{W}(v, x_{2})w} \nn
  &&= \res_{x_{0}}\res_{x_{2}}x_{0}^{p-l-i}x_{2}^{q+i}(x_{2}+x_{0})^{l}Y_{W}(Y(u, x_{0})v, x_{2}))
  w
\end{eqnarray}
for $i\ge k-q$.
Combining (\ref{defp}) and (\ref{upper-truncation-6}), we obtain
\begin{eqnarray}\label{eq1-3}
\lefteqn{\res_{x_{0}}\res_{x_{2}}\sum_{i=0}^{k-q-1}\binom{p-l}{i}
x_{0}^{p-l-i}x_{2}^{q+i} (x_{0}+x_{2})^{l}Y_W(u, x_{0}+x_{2})
Y_W(v, x_{2})w}\nn
&& = \res_{x_{0}}\res_{x_{2}}\sum_{i=0}^{k-q-1}\binom{p-l}{i}
x_{0}^{p-l-i}x_{2}^{q+i} (x_{2}+x_{0})^{l}Y_{W}(Y(u, x_{0})v, x_{2})w.
\end{eqnarray}

We now use induction on $k-q-1$ to prove
\begin{eqnarray}\label{eq1-3-0}
\lefteqn{\res_{x_{0}}\res_{x_{2}}x_{0}^{p-l}x_{2}^{q}(x_{0}+x_{2})^{l}Y_{W}(u, x_{0}+x_{2})
Y_{W}(v, x_{2})w} \nn
&&= \res_{x_{0}}\res_{x_{2}}x_{0}^{p-l}x_{2}^{q}(x_{2}+x_{0})^{l}Y_{W}(Y(u, x_{0})v, x_{2}))
w
\end{eqnarray}
for $p\in \Z$ and $q<k$.
When $k-q-1=0$, (\ref{eq1-3}) becomes
\begin{eqnarray*}
\lefteqn{\res_{x_{0}}\res_{x_{2}}x_{0}^{p-l}x_{2}^{q}(x_{0}+x_{2})^{l}Y_{W}(u, x_{0}+x_{2})
Y_{W}(v, x_{2})w} \nn
&&= \res_{x_{0}}\res_{x_{2}}x_{0}^{p-l}x_{2}^{q}(x_{2}+x_{0})^{l}Y_{W}(Y(u, x_{0})v, x_{2}))
w.
\end{eqnarray*}
Assume that (\ref{eq1-3-0}) holds when $0\le k-q-1<n$. When
$k-q-1=n$, $0\le k-q-i-1<n$ for $i=1, \dots, n=k-q-1$. Since $p$ is arbitrary,
we can replace $p$ by $p-i$ for any $i\in \Z$ in (\ref{eq1-3-0}).
Thus by the induction assumption,
\begin{eqnarray}\label{eq1-3-1}
\lefteqn{\res_{x_{0}}\res_{x_{2}}
x_{0}^{p-l-i}x_{2}^{q+i} (x_{0}+x_{2})^{l}Y_W(u, x_{0}+x_{2})
Y_W(v, x_{2})w}\nn
&& = \res_{x_{0}}\res_{x_{2}}
x_{0}^{p-l-i}x_{2}^{q+i} (x_{2}+x_{0})^{l}Y_{W}(Y(u, x_{0})v, x_{2})w.
\end{eqnarray}
for $i=1, \dots, n=k-q-1$.
From (\ref{eq1-3-1}) for $i=1, \dots, n=k-q-1$ and (\ref{eq1-3}),
we obtain
\begin{eqnarray*}
\lefteqn{\res_{x_{0}}\res_{x_{2}}
x_{0}^{p-l}x_{2}^{q}((x_{0}+x_{2})^{l}(Y_W(u, x_{0}+x_{2})
Y_W(v, x_{2}))w}\nn
&&=\res_{x_{0}}\res_{x_{2}}\sum_{i=0}^{k-q-1}\binom{p-l}{ i}
x_{0}^{p-l-i}x_{2}^{q+i}(x_{0}+x_{2})^{l}Y_W(u, x_{0}+x_{2})
Y_W(v, x_{2})w\nn
&&\;\;\;-\res_{x_{0}}\res_{x_{2}}\sum_{i=1}^{k-q-1}\binom{p-l}{ i}
x_{0}^{p-l-i}x_{2}^{q+i}(x_{0}+x_{2})^{l}Y_W(u, x_{0}+x_{2})
Y_W(v, x_{2})w\nn
&&=\res_{x_{0}}\res_{x_{2}}\sum_{i=0}^{k-q-1}\binom{p-l}{i}
x_{0}^{p-l-i}x_{2}^{q+i} (x_{2}+x_{0})^{l}Y_{W}(Y(u, x_{0})v, x_{2})w\nn
&&\;\;\;-\res_{x_{0}}\res_{x_{2}}\sum_{i=1}^{k-q-1}\binom{p-l}{i}
x_{0}^{p-l-i}x_{2}^{q+i} (x_{2}+x_{0})^{l}Y_{W}(Y(u, x_{0})v, x_{2})w\nn
&&= \res_{x_{0}}\res_{x_{2}}x_{0}^{p-l}x_{2}^{q}(x_{2}+x_{0})^{l}Y_{W}(Y(u, x_{0})v, x_{2})
w
\end{eqnarray*}
proving (\ref{eq1-3-0}) in this case.

Taking $i=0$ in  (\ref{upper-truncation-3}), we see that
(\ref{eq1-3-0}) also holds for $p\in \Z$ and $q\ge k$. Thus
(\ref{eq1-3-0}) holds for $p, q\in \Z$. But this means that
$$(x_{0}+x_{2})^{l}Y_W(u, x_{0}+x_{2})
Y_W(v, x_{2})w = (x_{2} + x_{0})^lY_{W}(Y(u, x_{0})v, x_{2})w,$$
that is, (\ref{wk-assoc}) holds.
\epfv

The following theorem from \cite{LL} says that the Jacobi identity for the module
follows from weak associativity for the module and skew symmetry for the vertex operator algebra:

\begin{thm}[\cite{LL}]\label{LL}
  Let $(V, Y, {\bf 1})$ be a triple that satisfies all the axioms in
  the definition of the notion of vertex algebra, in particular the
  skew symmetry property. Let $W$ be a vector space and let
  $Y_W(\cdot, x)$ be a linear map from $V$ to $(\edo W)[[x, x^{-1}]]$
  such that $Y_W({\bf 1}, x) = 1$ and $Y_W(v, x)w \in W((x))$ for $v
  \in V$ and $w \in W$. Assume that weak associativity holds for
  any $u, v \in V$ and $w \in W$, in the sense that there exists $l
  \in \N$ (depending on $u$ and $w$) such that
\begin{equation}
(x_{0}+x_{2})^{l}Y_W(u, x_0 + x_2)Y_W(v, x_2)w
=(x_{2}+x_{0})^{l}Y_W(Y(u, x_0)v, x_2)w.
\end{equation}
Then the Jacobi identity holds for $u, v \in V$ and $w \in W$.
\end{thm}

Therefore, by providing the two properties equivalent to the weak associativity for the module, given in Theorem
\ref{main theorem}, as well as the truncation
properties and skew symmetry for the vertex operator algebra, we can
obtain the major axiom --Jacobi identity for the vertex operator
algebra modules:

\begin{thm}\label{main thm}
Let $V$ be a $\Z$-graded vertex algebra, $W=\coprod_{n\in \N}W_{(n)}$
an $\N$-graded vector space and
$Y_W$ a linear map from $V\otimes W$ to $W((x))$. For $v\in V$ and $w\in W$,
we denote the image of $u\otimes w$ under $Y_W$
by $Y_W(u, x)w$ and $\res_{x}x^{n}Y_W(u, x)w$ by $u_{n}w$.
Assume that $u_{n}$ maps $W_{(k)}$ to
$W_{(k+m-n-1)}$ for $u\in V_{(m)}$ and $n\in \Z$ and $Y_{W}(\one, x)=1_{W}$.
Also assume that for $u, v\in V$, $w\in W$, there exist
$k, l \in \Z$ such that
(\ref{eq6}) and (\ref{eq7}) hold, and for $p, q\in \Z$, $u, v\in V$, $w\in W$,
the formulas (\ref{eq1-1}) and (\ref{upper-truncation-3}) hold.
Then $(W, Y_W)$ is an admissible $V$-module.\epf
\end{thm}

By a similar proof as Theorem \ref{LL}, we can also show the twisted
analogue of Theorem \ref{LL}:
\begin{thm}
  Let $(V, Y, {\bf 1})$ be $\Z$-graded vertex algebra. Suppose that
  $V$ has an automorphism $g$ of order $T$ and $V$ has an eigenspace
  decomposition with respect to the action of $g$ as
\[
V = \coprod_{r = 0}^{T-1} V^r,
\]
where
\[
V^r = \{v \in V | gv = e^{2\pi ir/T}v\}.
\]
Let $W$ be a vector space and let $Y_W(\cdot, x)$ be a linear map from
$V$ to $(\edo W)[[x^{\frac{1}{T}}, x^{-\frac{1}{T}}]]$ such that $Y_W({\bf 1}, x) = 1$ and
$Y_W(v, x)w \in W((x^{\frac{1}{T}}))$ for $v \in V$ and $w \in W$. Assume that
weak associativity holds for any $u \in V^{r}$, $v \in V$ and $w \in
W$, in the sense that there exists $l \in \frac{r}{T}+\N$ (depending on $u$ and
$w$) such that
\begin{equation}
(x_{0}+x_{2})^{l}Y_W(u, x_0 + x_2)Y_W(v, x_2)w
=(x_{2}+x_{0})^{l}Y_W(Y(u, x_0)v, x_2)w.
\end{equation}
Then the twisted Jacobi identity holds for $u, v \in V$ and $w \in W$.
\end{thm}

\begin{thm}\label{main thm 1}
  Let $V$ be a $\Z$-graded vertex algebra with a finite order
  automorphism $g$, $W=\coprod_{n\in \frac{1}{T}\N}W_{(n)}$ be a $\frac{1}{T}\N$-graded vector
  space and $Y_W$ be a linear map from $V\otimes W$ to $W((x))$. For
  $u \in V^r$ and $w\in W$, we denote the image of $u\otimes w$ under
  $Y_W$ by $Y_W(u, x)w$ and $\res_{x}x^{n}Y_W(u, x)w$ by $u_{n}w$.
  Assume that $u_{n}$ maps $W_{(k)}$ to $W_{(k+m-n-1)}$ for $u\in
  V_{(m)}$ and $n\in \frac{r}{T} + \Z$ and $Y_{W}(\one, x)=1_{W}$.  Also assume that
  for $u, v\in V$, $w\in W$, there exist $k, l \in \frac{1}{T}\Z$ such that
  (\ref{eq6}) and (\ref{eq7}) hold, and for $p \in l + \Z$, $q \in k + \Z$, $u, v\in V$,
  $w\in W$, the formulas (\ref{eq1-1}) and (\ref{upper-truncation-3})
  hold.  Then $(W, Y_W)$ is an admissible $g$-twisted $V$-module.\epf
\end{thm}

The component form of (\ref{upper-truncation-3}) is
\begin{equation}\label{upper-truncation-4}
\sum_{j=0}^{\infty}{l\choose j}(u_{j+p-l-i}v)_{q-j+l+i}w=0
\end{equation}
for $i\ge k-q$. Set $N = p + q + 2$ and $m = p - l - i$,
(\ref{upper-truncation-4}) becomes
\begin{equation}\label{eq1-2}
\sum_{j=0}^{\infty}{l\choose j}(u_{j+m}v)_{N - j - m -2}w=0
\end{equation}
for $m \leq N - k - l - 2$.
We set $o(a) = a(\wt a - 1)$ for
homogeneous $a \in V$. Then by taking $N = \wt u + \wt v$ in
(\ref{eq1-2}), we have
\begin{equation}\label{general eq}
o\bigg(\sum_{j=0}^{\infty}{l\choose j}u_{j+m}v\bigg)w = 0,
\end{equation}
where $m \leq \wt u + \wt v - k - l - 2$.

\begin{cor}\label{associate}
Let $W$ be an admissible $V$-module, $u, v \in V$ and $w \in W_{(n)}$ for $n \in \N$. Then
\[
o\bigg(\sum_{j=0}^{\wt u + n}{\wt u + n \choose j}u_{j+m}v\bigg)w = 0,
\]
where $m \leq - 2n-2$.
\end{cor}
\pf In equation (\ref{general eq}), let k = $\wt v + n$ and $l = \wt u
+ n$. \epfv

\begin{cor}\label{twisted}
  Suppose that $V$ has an automorphism $g$ of order $T$ and an
  eigenspace decomposition with respect to the action of $g$ as
\[
V = \coprod_{r =0}^{T-1} V^r,
\]
where
\[
V^r = \{v \in V | gv = e^{2\pi ir/T}v\}.
\]
Let $W$ be an admissible $g$-twisted $V$-module, $u \in V^{r}, v \in
V$ and $w \in W_{(0)}$. Then
\[
o\bigg(\sum_{j=0}^{\infty}{\wt u -1 +
  \delta_r + \frac{r}{T}\choose j}u_{j+m}v\bigg)w = 0
\]
for $m \leq -\delta_r - 1$.
\end{cor}
\pf In equation (\ref{general eq}), let k = $\wt v - \frac{r}{T}$ and
$l = \wt u -1 + \delta_r + \frac{r}{T}$. \epfv

The product formula (\ref{eq1-1}) will be used to calculate the product
of operator components acting on the module:
\begin{lemma}\label{lemma1}
In the setting of Theorem \ref{main theorem}. The product formula
  (\ref{eq1-1}) is equivalent to
\begin{eqnarray}\label{eq1sec3}
u_pv_qw = \res_{x_0}\res_{x_2}\left(\sum_{i =0}^{k-q-1}{p-l\choose i}x_{0}^{p-l-i}x_{2}^{i}\right)x_2^q(x_2+x_0)^lY_W(Y(u, x_0)v, x_2)w.\nn
\end{eqnarray}
\end{lemma}

\setcounter{equation}{0}
\section{Zhu algebra}
In this section, we will recall the definition and some properties of
the Zhu algebra introduced in \cite{Z}, a generalization of the Zhu algebra
$A_n(V)$ defined in \cite{DLM1} for $n \in \N$ and twisted
generalization of the Zhu algebra $A_g(V)$ defined in \cite{DLM2} for a
finite order automorphism $g$ of $V$.

We first recall the definition of the Zhu algebra $A(V)$ for a vertex
operator algebra $V$.

\begin{defn}[\cite{Z}]{\rm
Let $O(V)$ be the subspace of $V$ spanned by elements of the form
\[
\res_x \frac{(1 + x)^{\wt u}Y(u, x)v}{x^2}
\]
for homogeneous $u, v \in V$. Zhu algebra $A(V)$ is defined to be the
quotient space $V/O(V)$.  }
\end{defn}


The Zhu algebra was then generalized to an associative algebra $A_n(V)$
for $n \in \N$ and to a twisted Zhu algebra $A_g(V)$ for
an automorphism $g$ of $V$.

\begin{defn}[\cite{DLM1}]{\rm
For $n \in \N$, let $O_n(V)$ be the subspace of $V$ spanned by elements of the form
\[
\res_x \frac{(1 + x)^{\wt u+n}Y(u, x)v}{x^{2n+2}}
\]
for homogeneous $u, v \in V$. $A_n(V)$ is defined to be the quotient
space $V/O_n(V)$.}
\end{defn}

The subspace $O_n(V)$ is a two-sided ideal of $V$
and $A_n(V)$ is an associative algebra under the multiplication $*$
defined by
\[
u *_n v = \sum_{m=0}^n(-1)^m\binom{m+n}{n}\res_x \frac{(1 + x)^{\wt u
    + n}Y(u, x)v}{x^{n+m+1}},
\]
for homogeneous $u, v \in V$, and for general $u, v \in V$, $*$ is
defined by linearity. Also, for every homogeneous element $u \in V$
and $m \geq k \geq 0$, elements of the form
\begin{equation}
\res_x \frac{(1 + x)^{\wt u + n+ k}Y(u, x)v}{x^{m+2n+2}}
\end{equation}
lie in $O_n(V)$.

Let $W$ be an admissible $V$-module and let $\Omega_m(W)$ denote the subspace consisting of $m$th lowest
weight vectors in $W$, that is
\[
\Omega_m(W) = \{w \in W | u_nw = 0\;\;\;{\rm if}\; \wt u_n < -m\}.
\]
It was shown in \cite{DLM1} that $\Omega_n(W)$ is an $A_n(V)$-module
via the action $o(v+O_n(V)) = v_{\wt v -1}$ for $v \in V$. The first
necessary condition to prove is that $o(O_n(V))$ annihilates $\Omega_n(W)$. Corollary
\ref{associate} gives an alternating proof of this condition and also
motivates us to give another construction of functors between
categories of $A_n(V)$-modules and categories of admissible
$V$-modules (see Section $4$ and $5$ for the detail).

\begin{defn}[\cite{DLM2}]{\rm
For an automorphism $g$ of $V$ of finite order $T$, let $O_g(V)$ be the subspaceof $V$ spanned by elements of the form
\[
\res_x \frac{(1 + x)^{\wt u-1 + \delta_r + \frac{r}{T}}Y(u, x)v}{x^{1+\delta_r}}
\]
for homogeneous $u \in V^{r}, v \in V$, where
\[
V^r = \{v \in V | gv = e^{2\pi ir/T}v\}.
\]
The twisted Zhu algebra $A_g(V)$ is defined to be the quotient space
$V/O_g(V)$.}
\end{defn}

The subspace $O_g(V)$ is a two-sided ideal of $V$
and $A_g(V)$ is an associative algebra under the multiplication $*$
defined by
\begin{eqnarray*}
u *_g v =
\begin{cases}
\res_x \frac{(1 + x)^{\wt u}Y(u, x)v}{x}\;\;\;&\mbox{if}\; r =0,\\
0\;\;\;&\mbox{if}\; r>0.
\end{cases}
\end{eqnarray*}
for homogeneous $u \in V^r, v \in V$, and for general $u, v \in V$,
$*_g$ is defined by linearity. Also, for every homogeneous element $u
\in V^r$ and $m \geq k \geq 0$, elements of the form
\begin{equation}\label{z}
\res_x \frac{(1 + x)^{\wt u -1 + \delta_r+ \frac{r}{T}+ k}Y(u, x)v}{x^{m+\delta_r+1}}
\end{equation}
lie in $O_g(V)$.

The subspace $V^r \subset O_g(V)$ for $r \neq 0$
and $A_g(V)$ is a quotient of $A(V^0)$. Therefore, an $A_g(V)$-module
can be lift to an $A(V^0)$-module.

Let $W$ be an admissible twisted $V$-module and let $\Omega(W)$ denote the subspace consisting of lowest weight
vectors in $W$, that is
\[
\Omega(W) = \{w \in W| u_nw = 0\;\;\;{\rm if}\; \wt u_n < 0\}.
\]
It was shown in \cite{DLM2} that $\Omega(W)$ is an $A_g(V)$-module via
the action $o(v+O_g(V)) = v_{\wt v -1}$ for $v \in V^0$. The first
thing to show is that $o(O_g(V))$ annihilates $\Omega(W)$. Corollary
\ref{twisted} gives an alternating proof of this fact and also
motivates us to give another construction of functors between
categories of $A_g(V)$-modules and categories of admissible
$g$-twisted $V$-modules (see Section $6$ for the detail).

\setcounter{equation}{0}
\section{A functor $S_n$ from the category of $A_n(V)$-modules to the category of $V$-modules}
In the remaining part of this paper, we will assume the vertex operator algebra $V$ is $\N$-graded.

In this section, we will start from an $A_n(V)$-module $W$ and construct a vector space $S_n(W)$ with a linear map $Y_{S_n(W)}$ from $V \otimes S_n(W)$ to $S_n(W)((x))$. Then we use Theorem \ref{main theorem} to show that weak associativity holds for the pair $(S_n(W), Y_{S_n(W)})$ and hence the pair is an admissible $V$-module by Theorem \ref{main thm}.

Consider the affinization $V[t, t^{-1}] = V \otimes \C[t, t^{-1}]$ of
$V$ and the tensor algebra $T(V[t, t^{-1}])$ generated by $V[t,
t^{-1}]$. For simplicity, we shall denote $u \otimes t^m$ for $u \in
V$ and $m \in \Z$ by $u(m)$ and we shall omit the tensor product sign
$\otimes$ when we write an element of $T(V[t, t^{-1}])$. Thus $T(V[t,
t^{-1}])$ is spanned by elements of the form $u_1(m_1)\cdots u_k(m_k)$
for $u_i \in V$ and $m_i \in \Z$, $i = 1, \dots, k$.

Consider $T(V[t, t^{-1}])\otimes W$. Again for simplicity we shall
omit the tensor product sign. So $T(V[t, t^{-1}])\otimes W$ is
spanned by elements of the form $u_1(m_1)\cdots u_k(m_k)w$ for
$u_i \in V$, $m_i \in \Z$, $i = 1, \dots, k$ and $w \in W$ and
for any $u \in V$, $m \in Z$, $u(m)$ acts from the left on
$T(V[t, t^{-1}])\otimes W$. For homogeneous $u_i \in V$,
$m_i \in \Z$, $ i = 1, \dots, k$ and $w \in W$, we define the
degree of elements in $T(V[t, t^{-1}])\otimes W$ as follows:
\[
\deg\; u_1(m_1) \cdots u_k(m_k)w = (\wt u_1 - m_1 - 1) + \cdots (\wt u_k - m_k - 1) + n.
\]
For any $u \in V$, let
\[
Y_t(u, x): T(V[t, t^{-1}]) \otimes W \longrightarrow
T(V[t, t^{-1}]) \otimes W[[x, x^{-1}]]
\]
be defined by
\[
Y_t(u, x) = \sum_{m \in \Z}u(m)x^{-m-1}.
\]
For a homogeneous element $u \in V$, let \[o_t(u) = u(\wt u - 1).\]
Using linearity, we extend $o_t(u)$ to non-homogeneous $u$.

Let $\rho: A_n(V) \rightarrow \edo W$ be a representation of
associative algebra $A_n(V)$. Let $\mathcal{I}$ be the $\Z$-graded
$T(V[t, t^{-1}])$-submodule of
$T(V[t, t^{-1}]) \otimes W$ generated by elements of the forms
$u(m)w$ ($u \in V$, $\wt u - m -1 + \deg w < 0$,
$w \in T(V[t, t^{-1}])\otimes W$), $o_t(u)w - \rho(u + O_n(V))w$
($u \in V$, $w \in W$) and
\begin{eqnarray}\label{eqn2sec3}
&&u(p)v(q)w - \sum_{i = 0}^{\wt v + \deg w +n-q -1}
\sum_{j = 0}^{\wt u + \deg w + n}\left(\begin{array}{c}
p - \wt u - \deg w - n \\ i\end{array}\right)\nn
&&\cdot \left(\begin{array}{c}
 \wt u + \deg w + n\\ j\end{array}\right)(u_{p-\wt u -\deg w-i+j-n}v)(q+\wt u + \deg w+i-j+n)w\nn
\end{eqnarray}
($u, v \in V$, $q \in \Z$ such that $\wt v - q - 1 + \deg w \geq 0$,
$w \in T(V[t, t^{-1}]) \otimes W$). Note that relation (\ref{eqn2sec3}) is derived from formula (\ref{eq1sec3}) (also (\ref{eq1-1})) when $k = \wt v + \deg w + n$ and $l = \wt u + \deg w + n$.

Let
$\tilde{S_n}(W) = T(V[t, t^{-1}]) \otimes W/ \mathcal{I}$. Then $\tilde{S_n}(W)$
is also a $\Z$-graded $T(V[t, t^{-1}])$-module. In fact, by
definition of $\mathcal{I}$, we see that $\tilde{S_n}(W)$ is spanned
by elements of the form $u(m)w + \mathcal{I}$ for homogeneous
$u \in V$, $m \in \Z$ such that $m < \wt u + n - 1$ and $w \in W$.
In particular, we see that $\tilde{S_n}(W)$ has an $\N$-grading. Note
that $\mathcal{I} \cap W = \{0\}$, $W$ can be embedded into
$\tilde{S_n}(W)$ and $(\tilde{S_n}(W))_n = W$.

Let $\mathcal{J}$ be the $\N$-graded $T(V[t, t^{-1}])$-submodule
of $\tilde{S_n}(W)$ generated by
\begin{equation}\label{newequation4}
\sum_{j = 0}^{\wt u + \deg w + n}\left(\begin{array}{c} \wt u + \deg w + n\\
j\end{array}\right)(u_{j+m}v)(N-j-m-2)w
\end{equation}
($u, v \in V$, $w \in \tilde{S_n}(W)$, $N \in \Z$, $m \leq N-2 - \wt u -
\wt v -2n - 2\deg w$). Note that relation (\ref{newequation4}) comes from the formula (\ref{eq1-2}) by specializing $k = \wt v + \deg w + n$ and $l = \wt u + \deg w + n$.

Let $S_n(W) = \tilde{S_n}(W)/\mathcal{J}$. Then $S_n(W)$
is also an $\N$-graded $T(V[t, t^{-1}])$-module. We can still
use elements of $T(V[t, t^{-1}])\otimes W$ to represent elements
of $S_n(W)$. But note that these elements now satisfy relations.
We equip $S_n(W)$ with the vertex operator map
\[
Y_{S_n(W)}: V \otimes S_n(W) \longrightarrow S_n(W)[[x, x^{-1}]]
\]
given by
\[
u \otimes w \rightarrow Y(u, x)w = Y_t(u, x)w.
\]

\begin{thm}
The pair
$(S_n(W), Y_{S_n(W)})$ is an admissible $V$-module.
\end{thm}
\pf As in $\tilde{S_n}(W)$, for $u \in V$ and $w \in S_n(W)$, we also
have $u(m)w = 0$ when $m > \wt u + \deg w - 1$. Clearly,
\[
Y({\bf 1}, x) = I_{S_n(W)},
\]
where $I_{S_n(W)}$ is the identity operator on $S_n(W)$.

Since $S_n(W)$ satisfies the product formula and iterate formula in
Theorem \ref{main theorem}, where we specialize $k = \wt v + \deg w + n$ and $l = \wt u
+ \deg w + n$, by definition of $\mathcal{I}$ and
$\mathcal{J}$, $S_n(W)$ is an admissible weak $V$-module. \epfv

Let $W_{1}$ and $W_{2}$ be an $A_n(V)$-modules and
$f: W_{1} \to W_{2}$ a module map. Then $f$ induces a linear map
from $T(V[t, t^{-1}])\otimes W_{1}$ to $T(V[t, t^{-1}])\otimes W_{2}$.
By definition, this induced linear map in turn induces a linear map
$S_n(f)$ from $S_n(W_{1})$ to $S_n(W_{2})$. Since $Y_{S_n(W_{1})}$ and $Y_{S_n(W_{2})}$
are induced by $Y_{t}$ on $T(V[t, t^{-1}]) \otimes W_{1}$ and
$T(V[t, t^{-1}]) \otimes W_{2}$, respectively, we have
$$S_n(f)(Y_{S_n(W_{1})}(u, x)w_{1})= Y_{S_n(W_{2})}(u, x)S_n(f)(w_{1})$$
for $u\in V$ and $w_{1}\in S_n(W_{1})$. Thus $S_n(f)$ is a module map.
The following result is now clear:

\begin{cor}\label{functor-2}
  Let $V$ be an $\N$-graded vertex algebra.  Then the correspondence
  sending an $A_n(V)$-module $W$ to an admissible $V$-module $(S_n(W),
  Y_{S_n(W)})$ and an $A_n(V)$-module map $W_{1}\to W_{2}$ to a
  $V$-module map $S_n(f): S_n(W_{1})\to S_n(W_{2})$ is a functor from
  the category of $A_n(V)$-modules to the category of admissible
  $V$-modules.
\end{cor}

\setcounter{equation}{0}
\section{A universal property for $S_n$}
In this section, we prove that $S_n$ satisfies a natural universal
property and thus is the same as the functor constructed in
\cite{DLM1}. In particular, we achieve our goal of constructing
admissible $V$-modules from $A_n(V)$-modules without dividing
relations corresponding to the commutator formula for weak modules.

We use the following two lemmas to prove that $\mathcal{J} \cap W = 0$
in $\tilde{S_n}(W)$ and hence $(S_n(W))_n = W$.

\begin{lemma}\label{lemma2}
Let $M \in \N$ and $M \leq n$. Then in $\tilde{S_n}(W)$,
\[
a(\wt a -M-1)\bigg( \sum_{j=0}^{\wt u + n}\binom{\wt u +
  n}{j}(u_{j+m}v)(\wt u_{j+m}v + M -1)\bigg)w = 0,
\]
where $w \in W$ and $m \leq M-3n -2$.
\end{lemma}
\pf We apply formula (\ref{eq1sec3}) to
\[
a(\wt a - M - 1)\left(\sum_{j = 0}^{\wt u + n}\left(\begin{array}{c}
      \wt u + n \\ j\end{array}\right)(u_{j+m}v)(\wt u + \wt v - j -
  m + M -2)\right)w
\]
by specifying
\begin{eqnarray*}
&&p = \wt a - M - 1\\
&&q = \wt u + \wt v - j - m + M -2\\
&&k = \wt u + \wt v - j - m -1 + 2n\\
&&l = \wt a + 2n,
\end{eqnarray*}
we have
\begin{eqnarray*}
&&a(\wt a - M - 1)\left(\sum_{j = 0}^{\wt u + n}\left(\begin{array}{c} \wt u + n\\j\end{array}\right)(u_{j+m}v)(\wt u + \wt v - j - m + M -2)\right)w\nn
 &=&\res_{x_0}\res_{x_2}
  \sum_{i = 0}^{2n-M}\sum_{j = 0}^{\wt u + n}\binom{-M-2n-1}{i}\binom{\wt u + n}{j}x_0^{-M-2n-i-1}x_2^{\wt u + \wt v - j - m + M -2+i}\nn
  &&\;\;\;\;\;\;\;\cdot (x_2+x_0)^{\wt a + 2n}Y_{\tilde{S}_n(W)}(Y(a, x_0)(u_{j+m}v), x_2)w\nn
  &=& \res_y\res_{x_0}\res_{x_2}\sum_{i = 0}^{2n-M}\binom{-M-2n-1}{i}(1 + y)^{\wt u+n}y^mx_0^{-M-2n-i-1}\nn
  &&\;\;\;\;\;\;\;\cdot (x_2+x_0)^{\wt a +2n}Y_{\tilde{S}_n(W)}(Y(a, x_0)x_2^{L(0)+M-1+i}Y(u, y)v, x_2)w\nn
  &=&\res_y\res_{x_0}\res_{x_2}\sum_{i = 0}^{2n-M}\binom{-M-2n-1}{i}(1 + y)^{\wt u + n}y^mx_0^{-M-2n-i-1}x_2^{\wt u+\wt v+M-1+i}\nn
  &&\;\;\;\;\;\;\;\cdot (x_2+x_0)^{\wt a + 2n}Y_{\tilde{S}_n(W)}(Y(a, x_0)Y(u, x_2y)v, x_2)w\nn
   &=&\res_y\res_{x_0}\res_{x_2}\sum_{i = 0}^{2n-M}\binom{-M-2n-1}{i}(1 + y)^{\wt u + n}y^mx_0^{-M-2n-i-1}x_2^{\wt u+\wt v+M-1+i}\nn
   &&\;\;\;\;\;\;\;\cdot (x_2+x_0)^{\wt a + 2n}Y_{\tilde{S}_n(W)}(Y(u, x_2y)Y(a, x_0)v, x_2)w\nn
   &&+ \;\res_y\res_{x_0}\res_{x_2}\res_{x_1}\sum_{i = 0}^{2n-M}\binom{-M-2n-1}{i}(1 + y)^{\wt u + n}y^mx_0^{-M-2n-i-1}\nn
   &&\;\;\;\;\cdot x_2^{\wt u+\wt v+M-1+i}(x_2+x_0)^{\wt a + 2n} x_0^{-1}\delta(\frac{x_2y+x_1}{x_0})Y_{\tilde{S}_n(W)}(Y(Y(a, x_1)u, x_2y)v, x_2)w.\nn
 \end{eqnarray*}

By examining the monomials in $y$ in the first term of the right-hand side, we know that the first term of the right-hand side is a sum of elements of form
\[
o\left(\sum_{j = 0}^{\wt u + n}{\wt u+n \choose j}(u_{j + m}\tilde{a})\right) w
\]
with $m \leq M-3n-2$ for some $\tilde{a} \in V$. Hence the first term
is an action of elements in $O_n(V)$ on $w$ which is $0$. We only need
to prove the second term is also an action of a sum of elements of
form $O_n(V)$ on $w$. The second term equals
\begin{eqnarray*}
&&\res_y\res_{x_0}\res_{x_2}\res_{x_1}\sum_{i = 0}^{2n-M}\binom{-M-2n-1}{i}(1 + y)^{\wt u + n}y^mx_0^{-M-2n-i-1}\nn
   &&\;\;\;\;\cdot x_2^{\wt u+\wt v+M-1+i}(x_2+x_0)^{\wt a + 2n} x_0^{-1}\delta(\frac{x_2y+x_1}{x_0})Y_{\tilde{S}_n(W)}(Y(Y(a, x_1)u, x_2y)v, x_2)w\nn
&=&\res_y\res_{x_2}\res_{x_1}\sum_{i = 0}^{2n-M}\binom{-M-2n-1}{i}(1 + y)^{\wt u + n}y^m(x_2y+x_1)^{-M-2n-i-1}\nn
   &&\;\;\;\;\cdot x_2^{\wt u+\wt v+M-1+i}(x_2+x_2y+x_1)^{\wt a + 2n} Y_{\tilde{S}_n(W)}(Y(Y(a, x_1)u, x_2y)v, x_2)w\nn
&=&\res_y\res_{x_2}\res_{x_3}\sum_{i = 0}^{2n-M}\binom{-M-2n-1}{i}(1 + y)^{\wt u + \wt a+3n+1}y^m(x_2y+x_3(1+y))^{-M-2n-i-1}\nn
   &&\;\;\;\;\cdot x_2^{\wt u+\wt v+M-1+i}(x_2+x_3)^{\wt a + 2n} Y_{\tilde{S}_n(W)}(Y(Y(a, (1+y)x_3)u, x_2y)v, x_2)w\nn
&=&\res_y\res_{x_2}\res_{x_3}\sum_{i = 0}^{2n-M}\sum_{j = 0}^{\infty}\binom{-M-2n-1}{i}\binom{-M-2n-i-1}{j}(1 + y)^{\wt u + \wt a+3n+1+j}\nn
   &&\;\;\;\;\cdot y^{m-M-2n-i-j-1}x_2^{\wt u+\wt v-2n-j-2}x_3^j(x_2+x_3)^{\wt a + 2n} Y_{\tilde{S}_n(W)}(Y(Y(a, (1+y)x_3)u, x_2y)v, x_2)w\nn
&=&\res_y\res_{x_2}\res_{x_3}\sum_{i = 0}^{2n-M}\sum_{j = 0}^{\infty}\binom{-M-2n-1}{i}\binom{-M-2n-i-1}{j}x_2^{\wt u+\wt v-2n-j-2}x_3^j\nn
   &&\;\;\;\;\cdot (x_2+x_3)^{\wt a + 2n} Y_{\tilde{S}_n(W)}(y^{m-M-2n-i-j-1}Y((1 + y)^{L(0)+3n+1+j}Y(a, x_3)u, x_2y)v, x_2)w.\nn
\end{eqnarray*}
Since $m \leq M-3n-2$ and $i,j\geq 0$, we obtain that it is an
action of elements of the form
\[
\res_y y^{m'}Y(1+y)^{L(0)+3n+1+j}(Y(a, x_3)u, y)v
\]
on $w$, where $m'=m-M-2n-i-j-1 \leq -5n-i-j-3$. Apparently, this is an element in $O_n(V)$, the action equals $0$. \epfv

We proceed to prove the next proposition:
\begin{prop}\label{lemma3}
  Let $M, N \in \Z$ such that $0 \leq N+n$ and $M \leq N+n$. Then in
  $\tilde{S_n}(W)$,
\begin{eqnarray}\label{eqn3sec3}
&&a(\wt a - M + N - 1)\nn
&&\cdot\left(\sum_{j = 0}^{\wt u + N + 2n}\left(\begin{array}{c} \wt u + N + 2n\\ j\end{array}\right)(u_{j+m}v)(\wt u + \wt v - j - m + M -2)\right)\nn
&&\cdot b(\wt b - N - 1)w = 0
\end{eqnarray}
where $m \leq M -2N -4n -2$ and $w \in W$.
\end{prop}
\pf By Lemma \ref{lemma2}, it suffices to show that the element
\[
\left(\sum_{j = 0}^{\wt u + N + 2n}\left(\begin{array}{c} \wt u + N +
      2n\\ j\end{array}\right)(u_{j+m}v)(\wt u + \wt v - j - m + M
  -2)\right)\cdot b(\wt b - N - 1)w
\]
for $m \leq M -2N -4n -2$, is of the form
\[
\sum_{j = 0}^{\wt \tilde{u} + n}\left(\begin{array}{c} \wt \tilde{u}
    + n\\ j\end{array}\right)(\tilde{u}_{j+m}\tilde{v})(\wt \tilde{u}
+ \wt \tilde{v} - j - m + M - N -2)w
\]
for $m \leq M-N - 3n -2$, $\tilde{u}, \tilde{v} \in V$.

We apply formula (\ref{eq1sec3}) to
\[
\left(\sum_{j = 0}^{\wt u + N + 2n}\left(\begin{array}{c} \wt u + N +
      2n\\ j\end{array}\right)(u_{j+m}v)(\wt u + \wt v - j - m + M
  -2)\right) b(\wt b - N - 1)w
\]
by specifying
\begin{eqnarray*}
&&p = \wt u + \wt v - j - m + M -2\\
&&q = \wt b - N - 1\\
&&k = \wt b + 2n\\
&&l = \wt u + \wt v - j - m -1 + 2n,
\end{eqnarray*}
we have
\begin{eqnarray*}
&&\sum_{j = 0}^{\wt u + N + 2n}{\wt u + N + 2n \choose j}(u_{j + m}v)(\wt u_{j+m}v + M - 1)b(\wt b -N-1)w\nn
&=& \res_{x_0}\res_{x_2}\res_y \sum_{i =0}^{2n+N}\binom{M-2n-1}{i}x_0^{M-2n-i-1}x_2^{\wt b-N-1+i}(1+y)^{\wt u + N + 2n}y^m\nn
&&\;\;\;\;\;\;\;\;\;\;\;\cdot Y_{\tilde{S}_n(W)}(Y((x_2+x_0)^{L(0)+2n}Y(u,y)v, x_0)b, x_2)w\nn
&=& \res_{x_0}\res_{x_2}\res_y \sum_{i =0}^{2n+N}\binom{M-2n-1}{i}x_0^{M-2n-i-1}x_2^{\wt b-N-1+i}(1+y)^{\wt u + N + 2n}y^m\nn
&&\;\;\;\;\;\;\;\;\;\;\;\cdot (x_2+x_0)^{\wt u+\wt v+2n}Y_{\tilde{S}_n(W)}(Y(Y(u,(x_2+x_0)y)v, x_0)b, x_2)w\nn
&=& \res_{x_0}\res_{x_2}\res_{x_3} \sum_{i =0}^{2n+N}\binom{M-2n-1}{i}x_0^{M-2n-i-1}x_2^{\wt b-N-1+i}(x_2+x_0+x_3)^{\wt u + N + 2n}\nn
&&\;\;\;\;\;\;\;\;\;\;\;\cdot x_3^m(x_2+x_0)^{\wt v-N-m-1}Y_{\tilde{S}_n(W)}(Y(Y(u,x_3)v, x_0)b, x_2)w\nn
&=& \res_{x_0}\res_{x_2}\res_{x_1} \sum_{i =0}^{2n+N}\binom{M-2n-1}{i}x_0^{M-2n-i-1}x_2^{\wt b-N-1+i}(x_2+x_1)^{\wt u + N + 2n}\nn
&&\;\;\;\;\;\;\;\;\;\;\;\cdot (x_1-x_0)^m(x_2+x_0)^{\wt v-N-m-1}Y_{\tilde{S}_n(W)}(Y(u,x_1)Y(v,x_0)b, x_2)w\nn
&& -\; \res_{x_0}\res_{x_2}\res_{x_1} \sum_{i =0}^{2n+N}\binom{M-2n-1}{i}x_0^{M-2n-i-1}x_2^{\wt b-N-1+i}(x_2+x_1)^{\wt u + N + 2n}\nn
&&\;\;\;\;\;\;\;\;\;\;\;\cdot (-x_0+x_1)^m(x_2+x_0)^{\wt v-N-m-1}Y_{\tilde{S}_n(W)}(Y(v, x_0)Y(u, x_1)b, x_2)w\nn
\end{eqnarray*}
By examining the monomials in $x_0$ in the second term of the
right-hand side, it is a sum of elements of the form
\[
\sum_{j = 0}^{\wt v -N - m -1}\left(\begin{array}{c} \wt v -N-m-1\\ j\end{array}\right)(v_{j+m'}\tilde{u})(\wt v + \wt
\tilde{u} - j - m + M - N -2)
\]
acting on $w$ for some $\tilde{u} \in V$ and $m' \leq m +M-2n-1$. Since
\[
-N - m -1 \geq -M+N+4n+1 \geq n,
\]
this element is of the form
\[
\sum_{j = 0}^{\wt v + n}\left(\begin{array}{c} \wt v + n \\
    j\end{array}\right)(v_{j+m''}\tilde{u})(\wt v + \wt \tilde{u} - j -
m + M - N -2)
\]
for $m'' =  m' - N - m - n - 1 \leq M-N-3n-2$ and $\tilde{u} \in
V$. It is easy to see by examining the monimials in $x_1$ in first
term of the right-hand side that the first term is also a sum of
elements of this form. \epfv

The following theorem is an easy consequence of Proposition \ref{lemma3}:
\begin{thm}\label{lemma2sec3}
In $\tilde{S_n}(W)$,
\[
\mathcal{J} \cap W = 0.
\]
The embedding of $W$ to $\tilde{S_n}(W)$ induces an injection $e_W$ of
$A_n(V)$-modules from $W$ to $\Omega_n(S_n(W))$.
\end{thm}

\begin{thm}
The functor $S_n$ has the following universal property:
Let $W$ be an $A_n(V)$-module. For any admissible $V$-module
$\tilde{W}$ and any
$A_n(V)$-module map $f: W \to \Omega_n(\tilde{W})$, there exists a unique
$V$-module map $\tilde{f}: S_n(W)\to \tilde{W}$ such that
$\tilde{f}|_{e_W(W)}=f\circ e_{W}^{-1}$.
\end{thm}

In \cite{DLM1}, a functor, denoted by $\bar{M}_{n}$,
was constructed explicitly and was proved to satisfy the same universal property
above. The following result achieves our goal
of constructing this functor without dividing
relations corresponding to the commutator formula for weak modules:

\begin{cor}
The functor $S_n$ is equal to
the functor $\bar{M}_{n}$ constructed in \cite{DLM1}.
\end{cor}
\pf
This result follows immediately from the universal property.
\epfv

\setcounter{equation}{0}
\section{A functor $S_g$ from the category of $A_g(V)$-modules to the category of twisted $V$-module}
In this section, we assume the vertex operator algebra $V$ has an automorphism $g$ with finite order $T$ and an
eigenspace decomposition with respect to the action of $g$ as
\[
V = \coprod_{r =0}^{T-1} V^r,
\]
where
\[
V^r = \{v \in V | gv = e^{2\pi ir/T}v\}.
\]

We shall start with an $A_g(V)$-module $W$ and construct a vector space $S_g(W)$ with a linear map $Y_{S_g(W)}$ from $V \otimes S_g(W)$ to $S_g(W)((x^{\frac{1}{T}}))$. Then we show that the pair $(S_g(W), Y_{S_g(W)})$ is an admissible twisted $V$-module by verifying that weak associativity holds on the pair.

Consider the affinization $V^g[t, t^{-1}] = \coprod_{r=0}^{T-1}V^r
\otimes t^{r/T}\C[t, t^{-1}]$ of $V$ and the tensor algebra $T(V^g[t,
t^{-1}])$ generated by $V^g[t, t^{-1}]$. For simplicity, we shall
denote $u \otimes t^m$ for $u \in V^r$ and $m \in \frac{r}{T} + \Z$ by
$u(m)$ and we shall omit the tensor product sign $\otimes$ when we
write an element of $T(V^g[t, t^{-1}])$. Thus $T(V^g[t, t^{-1}])$ is
spanned by elements of the form $u_1(m_1)\cdots u_k(m_k)$ for $u_i \in
V^{r_i}$ and $m_i \in \frac{r_i}{T} + \Z$, $r_i = 0, 1, \dots, T-1$
for $i = 1, \dots, k$.

Consider $T(V^g[t, t^{-1}])\otimes W$. Again for simplicity we
omit the tensor product sign. So $T(V^g[t, t^{-1}])\otimes W$ is
spanned by elements of the form $u_1(m_1)\cdots u_k(m_k)w$ for $u_i
\in V^{r_i}$, $m_i \in \frac{r_i}{T} + \Z$, $i = 1, \dots, k$ and $w
\in W$ and for any $u \in V^r$, $m \in \frac{r}{T} + \Z$, $u(m)$ acts
from the left on $T(V^g[t, t^{-1}])\otimes W$. For homogeneous $u_i
\in V^{r_i}$, $m_i \in \frac{r_i}{T} + \Z$, $ i = 1, \dots, k$ and $w
\in W$, we define the degree of elements in $T(V^g[t, t^{-1}])\otimes W$ as
\[
\deg \; u_1(m_1) \cdots u_k(m_k)w = (\wt u_1 - m_1 - 1) + \cdots (\wt u_k - m_k - 1).
\]
For any $u \in V^{r}$, define
\[
Y_t(u, x): T(V^g[t, t^{-1}]) \otimes W \longrightarrow
T(V^g[t, t^{-1}]) \otimes W[[x^{\frac{1}{T}}, x^{-\frac{1}{T}}]]
\]
to be
\[
Y_t(u, x) = \sum_{m \in \frac{r}{T} + \Z}u(m)x^{-m-1}.
\]
For a homogeneous element $u \in V^0$, let \[o_t(u) = u(\wt u - 1).\]
Using linearity, we extend $o_t(u)$ to non-homogeneous $u \in V^0$.

Let $\rho: A_g(V) \rightarrow \edo W$ be a representation of
associative algebra $A_g(V)$. Let $\mathcal{I}$ be the $\Z$-graded
$T(V^g[t, t^{-1}])$-submodule of
$T(V^g[t, t^{-1}]) \otimes W$ generated by elements of the forms
$u(m)w$ ($u \in V$, $\wt u - m -1 + \deg w < 0$,
$w \in T(V^g[t, t^{-1}])\otimes W$), $o_t(u)w - \rho(u + O_g(V)\cap V^0)w$
($u \in V^0$, $w \in W$) and
\begin{eqnarray*}
  &&u(p)v(q)w - \sum_{i = 0}^{\wt v + \deg w-q -1}
  \sum_{j = 0}^{\infty}\binom{p - \wt u - \deg w -\delta_r - \frac{r}{T}}{i}
      \binom{\wt u + \deg w + \delta_r + \frac{r}{T}}{j}\nn
  &&\cdot (u_{p-\wt u -\deg w - \delta_r - \frac{r}{T}-i+j}v)(q+\wt u + \deg w+ \delta_r + \frac{r}{T}+i-j)w
\end{eqnarray*}
($u \in V^r, v \in V$, $q \in \Z$ such that $\wt v - q - 1 + \deg w \geq 0$,
$w \in T(V^g[t, t^{-1}]) \otimes W$).

Let
$\tilde{S}(W) = T(V^g[t, t^{-1}]) \otimes W/ \mathcal{I}$. Then $\tilde{S}(W)$
is also a $\frac{1}{T}\Z$-graded $T(V^g[t, t^{-1}])$-module. In fact, by
definition of $\mathcal{I}$, we see that $\tilde{S}(W)$ is spanned
by elements of the form $u(m)w + \mathcal{I}$ for homogeneous
$u \in V^r$, $m \in \frac{r}{T}+\Z$ such that $m < \wt u - 1$ and $w \in W$.
In particular, we see that $\tilde{S}(W)$ has an $\frac{1}{T}\N$-grading. Note
that $\mathcal{I} \cap W = \{0\}$, $W$ can be embedded into
$\tilde{S}(W)$ and $(\tilde{S}(W))_0 = W$.

Let $\mathcal{J}$ be the $\frac{1}{T}\N$-graded $T(V^g[t, t^{-1}])$-submodule
of $\tilde{S}(W)$ generated by
\[
\sum_{j = 0}^{\infty}\left(\begin{array}{c} \wt u + \deg w+ \delta_r + \frac{r}{T}\\
    j\end{array}\right)(u_{j+m}v)(N-j-m-2)w
\]
($u \in V^r, v \in V^s$, $w \in \tilde{S}(W)$, $N \in \frac{r+s}{T} +
\Z$, $m \leq N-2 - \wt u - \wt v - 2\deg w - \delta_r - \frac{r}{T}$).

\begin{thm}\label{lemma4}
Let $\mathcal{J}$ be the $\frac{1}{T}\N$-graded $T(V^g[t, t^{-1}])$-submodule
of $\tilde{S}(W)$ defined above. Then in $\tilde{S}(W)$,
\[
\mathcal{J} \cap W = 0.
\]
\end{thm}

\pf It suffices to prove the element of the form
\begin{eqnarray}\label{eqn3sec33}
&&a(\wt a - M + N - 1)\nn
&&\cdot \left(\sum_{j = 0}^{\infty}\left(\begin{array}{c} \wt u + N + \delta_r + \frac{r}{T}\\ j\end{array}\right)(u_{j+m}v)(\wt u + \wt v - j - m + M -2)\right)\nn
  &&\cdot b(\wt b - N - 1)w
\end{eqnarray}
equals $0$ in $\tilde{S}(W)$, where $a \in V^d$, $u \in V^{r}$, $v \in V^s$,
$b \in V^t$ ($d+r+s+t \equiv 0\; {\rm mod}\; T$), $w \in W$, $M \in
\frac{r+s}{T} + \Z, N \in \frac{T-t}{T} + \Z$,
\[
m \leq -N - 2 - (N-M)- \delta_r - \frac{r}{T}.
\]
We will prove the claim in the following cases:

(i)$N < 0$.  Since the element $b(\wt b - N - 1)w$ has negative
grading $N$, the expression (\ref{eqn3sec33}) lies in $\mathcal{I}$ and
is $0$ in $\tilde{S}(W)$;

(ii)$N \geq 0$ and $M > N$.
Applying formula (\ref{eq1sec3}) to the expression
\[
\left(\sum_{j = 0}^{\infty}\left(\begin{array}{c} \wt u + N + \delta_r +
      \frac{r}{T}\\ j\end{array}\right)(u_{j+m}v)(\wt u + \wt v - j -
  m + M -2)\right) b(\wt b - N - 1)w,
\]
it becomes an element with negative grading $N-M$, which lies in
$\mathcal{I}$ and hence equals $0$ in $\tilde{S}(W)$.

(iii)$N \geq 0$ and $N \geq M$. We will show this in the remaining of the proof.

It is easy to see that formula (\ref{eq1sec3}) holds for $\tilde{S}(W)$. We
will use formula (\ref{eq1sec3}) to simplify (\ref{eqn3sec33}) and show
that it is actually an element of $O_g(V)$ acting on $w \in W$ and
hence equals $0$ in $\tilde{S}(W)$. Note that
\begin{eqnarray*}
\res_y(1 + y)^{\wt u + N + \delta_r + \frac{r}{T}}y^mY(u, y)v = \sum_{j = 0}^{\infty}\left(\begin{array}{c} \wt u + N + \delta_r + \frac{r}{T}\\ j\end{array}\right)u_{j+m}v.
\end{eqnarray*}
For simplicity, we say an element is of form $O_N$ if it can be written as a linear combination of elements of the form:
\[
\sum_{j = 0}^{\infty}{\wt \tilde{u} + N + \delta_{\tilde{r}} +
  \frac{\tilde{r}}{T}\choose j}(\tilde{u}_{j + n}\tilde{v})(\wt
\tilde{u}_{j+n}\tilde{v} + N - 1)\tilde{w}
\]
with $n \leq -N-\delta_{\tilde{r}} - \frac{\tilde{r}}{T}-2$ for some
$\tilde{u} \in V^{\tilde{r}}$, $\tilde{v} \in V$.

Applying formula (\ref{eq1sec3}) to
\begin{eqnarray}
&&a(\wt a - M + N - 1)\nn
&&\cdot\left(\sum_{j = 0}^{\infty}\left(\begin{array}{c} \wt u + N + \delta_r + \frac{r}{T}\\ j\end{array}\right)(u_{j+m}v)(\wt u + \wt v - j - m + M -2)\right)\tilde{w},\nn
\end{eqnarray}
where $\tilde{w} = b(\wt b - N - 1)w \in \tilde{S}(W)$, by specializing
\begin{eqnarray*}
&&p = \wt a - M + N - 1\\
&&q = \wt u + \wt v - j - m + M -2\\
&&k = \wt u + \wt v - j - m + N -1\\
&&l = \wt a + N + \delta_d + \frac{d}{T},
\end{eqnarray*}
we have
\begin{eqnarray*}
&&a(\wt a - M + N - 1)\nn
&&\cdot\left(\sum_{j = 0}^{\infty}\left(\begin{array}{c} \wt u + N + \delta_r + \frac{r}{T}\\ j\end{array}\right)(u_{j+m}v)(\wt u + \wt v - j - m + M -2)\right)\tilde{w}\nn
& = & \sum_{j = 0}^{\infty}\left(\begin{array}{c} \wt u + N + \delta_r + \frac{r}{T}\\ j\end{array}\right)\res_{x_0}\res_{x_2}\sum_{i = 0}^{N-M}x_0^{-M-i-\delta_d - \frac{d}{T}-1}\nn
&&\;\;\;\;\;\;\;\cdot (x_2+x_0)^{\wt a + N + \delta_d + \frac{d}{T}}x_2^{\wt u + \wt v - j - m + M -2+i}Y_{\tilde{S}(W)}(Y(a, x_0)(u_{j+m}v), x_2)\tilde{w}\nn
&=& \res_y\res_{x_0}\res_{x_2}\sum_{i = 0}^{N-M}(1 + y)^{\wt u + N+ \delta_r + \frac{r}{T}}y^mx_0^{-M-i-\delta_d - \frac{d}{T}-1}\nn
&&\;\;\;\;\;\;\;\cdot (x_2+x_0)^{\wt a + N + \delta_d + \frac{d}{T}}Y_{\tilde{S}(W)}(Y(a, x_0)x_2^{L(0)+M-1+i}Y(u, y)v, x_2)\tilde{w}\nn
&=& \res_y\res_{x_0}\res_{x_2}\sum_{i = 0}^{N-M}(1 + y)^{\wt u + N + \delta_r + \frac{r}{T}}y^mx_0^{-M-i-\delta_d - \frac{d}{T}-1}x_2^{\wt u+\wt v+M-1+i}\nn
&&\;\;\;\;\;\;\;\cdot (x_2+x_0)^{\wt a + N + \delta_d + \frac{d}{T}}Y_{\tilde{S}(W)}(Y(u, x_2y)Y(a, x_0)v, x_2)\tilde{w}\nn
&& +\; \res_y\res_{x_0}\res_{x_2}\res_{x_1}\sum_{i = 0}^{N-M}(1 + y)^{\wt u + N + \delta_r + \frac{r}{T}}y^mx_0^{-M-i-\delta_d - \frac{d}{T}-1}x_2^{\wt u+\wt v+M-1+i}\nn
&&\;\;\;\;\;\;\;\cdot (x_2+x_0)^{\wt a + N + \delta_d + \frac{d}{T}}x_0^{-1}\delta(\frac{x_2y+x_1}{x_0})Y_{\tilde{S}(W)}(Y(Y(a, x_1)u, x_2y)v, x_2)\tilde{w}.
\end{eqnarray*}

By examining the monomials in $y$ in the first term of the right-hand
side, we know that the first term of the right-hand side is a sum of
elements of form $O_N$. We only need to prove the second term is also
a sum of elements of form $O_N$. The second term equals
\begin{eqnarray*}
&&\res_y\res_{x_0}\res_{x_2}\res_{x_1}\sum_{i = 0}^{N-M}(1 + y)^{\wt u + N + \delta_r + \frac{r}{T}}y^mx_0^{-M-i-\delta_d - \frac{d}{T}-1}x_2^{\wt u+\wt v+M-1+i}\nn
&&\;\;\;\;\;\;\;\cdot (x_2+x_0)^{\wt a + N + \delta_d + \frac{d}{T}}x_0^{-1}\delta(\frac{x_2y+x_1}{x_0})Y_{\tilde{S}(W)}(Y(Y(a, x_1)u, x_2y)v, x_2)\tilde{w}\nn
&=&\res_y\res_{x_2}\res_{x_1}\sum_{i = 0}^{N-M}(1 + y)^{\wt u + N + \delta_r + \frac{r}{T}}y^m(x_2y+x_1)^{-M-i-\delta_d - \frac{d}{T}-1}x_2^{\wt u+\wt v+M-1+i}\nn
&&\;\;\;\;\;\;\;\cdot (x_2+x_2y+x_1)^{\wt a + N + \delta_d + \frac{d}{T}}Y_{\tilde{S}(W)}(Y(Y(a, x_1)u, x_2y)v, x_2)\tilde{w}\nn
&=&\res_y\res_{x_2}\res_{x_3}\sum_{i = 0}^{N-M}(1 + y)^{\wt u + \wt a+2N + \delta_r+\delta_d + \frac{r+d}{T}+1}y^m(x_2y+x_3(1+y))^{-M-i-\delta_d - \frac{d}{T}-1}\nn
&&\;\;\;\;\;\;\;\cdot x_2^{\wt u+\wt v+M-1+i}(x_2+x_3)^{\wt a + N + \delta_d + \frac{d}{T}}Y_{\tilde{S}(W)}(Y(Y(a, (1+y)x_3)u, x_2y)v, x_2)\tilde{w}\nn
&=&\res_y\res_{x_2}\res_{x_3}\sum_{i = 0}^{N-M}\sum_{j=0}^{\infty}\binom{-M-i-\delta_d - \frac{d}{T}-1}{j}(1 + y)^{\wt u + \wt a+2N + \delta_r+\delta_d + \frac{r+d}{T}+1+j}\nn
&&\;\;\;\;\;\;\;\cdot x_3^jy^{m-M-i-\delta_d - \frac{d}{T}-1-j}x_2^{\wt u+\wt v-\delta_d - \frac{d}{T}-2-j}(x_2+x_3)^{\wt a + N + \delta_d + \frac{d}{T}}\nn
&&\;\;\;\;\;\;\;\cdot Y_{\tilde{S}(W)}(Y(Y(a, (1+y)x_3)u, x_2y)v, x_2)\tilde{w}\nn
&=&\res_y\res_{x_2}\res_{x_3}\sum_{i = 0}^{N-M}\sum_{j=0}^{\infty}\binom{-M-i-\delta_d - \frac{d}{T}-1}{j}\nn
&&\;\;\;\;\;\;\;\cdot x_3^j x_2^{\wt u+\wt v-\delta_d - \frac{d}{T}-2-j}(x_2+x_3)^{\wt a + N + \delta_d + \frac{d}{T}}\nn
&&\;\;\;\;\;\;\;\cdot Y_{\tilde{S}(W)}(y^{m-M-i-\delta_d - \frac{d}{T}-1-j}Y((1 + y)^{L(0)+2N + \delta_r+\delta_d + \frac{r+d}{T}+1+j}Y(a, x_3)u, x_2y)v, x_2)\tilde{w}.\nn
\end{eqnarray*}
By checking the monomial in $y$ of the right-hand side, it is the sum of elements of the form
\begin{equation}\label{newequation6}
\res_y y^{m'}Y(1+y)^{L(0)+2N + \delta_r+\delta_d + \frac{r+d}{T}+1+j}(Y(a, x_3)u, y)v,
\end{equation}
where $m' \leq m-M-i-\delta_d - \frac{d}{T}-1-j \leq -2N-\delta_r-\delta_d-\frac{r+d}{T}-3-j$. Note that
\[
2N + \delta_r+\delta_d + \frac{r+d}{T}+1+j \geq N + \delta_{\overline{r+d}}+ \frac{\overline{r+d}}{T},
\]
here we use $\overline{r+d}$ to denote the residue of $r+d$ modulo $T$, the element (\ref{newequation6}) can be written as an element of the form
\[
\res_y y^{m''}Y(1+y)^{L(0)+N + \delta_{\overline{r+d}}+ \frac{\overline{r+d}}{T}}(Y(a, x_3)u, y)v,
\]
where $m'' \leq m' + 2N + \delta_r+\delta_d + \frac{r+d}{T}+1+j - (N + \delta_{\overline{r+d}}+ \frac{\overline{r+d}}{T}) \leq -N-\delta_{\overline{r+d}}-\frac{\overline{r+d}}{T}-2$. Thus the element (\ref{newequation6}) is of form $O_N$. We proved that
\begin{eqnarray*}
&&a(\wt a - M + N - 1)\nn
&&\cdot\left(\sum_{j = 0}^{\infty}\left(\begin{array}{c} \wt u + N + \delta_r + \frac{r}{T}\\ j\end{array}\right)(u_{j+m}v)(\wt u + \wt v - j - m + M -2)\right)\tilde{w}
\end{eqnarray*}
is a sum of elements of the form
\[
\sum_{j = 0}^{\infty}{\wt \tilde{u} + N + \delta_{\tilde{r}} +
  \frac{\tilde{r}}{T}\choose j}(\tilde{u}_{j + n}\tilde{v})(\wt
\tilde{u}_{j+n}\tilde{v} + N - 1)\tilde{w}
\]
with $n \leq -N-\delta_{\tilde{r}} - \frac{\tilde{r}}{T}-2$ for some
$\tilde{u} \in V^{\tilde{r}}$, $\tilde{v} \in V$. For simplicity, we
still write this element as
\begin{equation}\label{eqn4sec3}
\sum_{j = 0}^{\infty}{\wt u + N + \delta_r + \frac{r}{T} \choose j}(u_{j + m}v)(\wt u_{j+m}v + N - 1)b(\wt b -N-1)w
\end{equation}
for $m \leq -N-\delta_r - \frac{r}{T}-2$, $u \in V^{r}$ and $w \in
W$. 

We shall prove that (\ref{eqn4sec3}) is $0$ in $\tilde{S}(W)$. Applying formula (\ref{eq1sec3}) to the expression (\ref{eqn4sec3}),
by specializing
\begin{eqnarray*}
&& p = \wt u_{j+m}v + N -1\nn
&& q = \wt b -N - 1\nn
&& k = \wt b\nn
&& l = \wt u_{j+m}v + \delta_{\overline{r+s}} + \frac{\overline{r+s}}{T},
\end{eqnarray*}
here we use $\overline{r+s}$ to denote the residue of $r+s$ modulo $T$, we have
\begin{eqnarray*}
&&\sum_{j = 0}^{\infty}{\wt u + N+ \delta_r + \frac{r}{T} \choose j}(u_{j + m}v)(\wt u_{j+m}v + N - 1)b(\wt b -N-1)w\nn
&=& \res_{x_0}\res_{x_2}\res_y \sum_{i =0}^Nx_0^{N-i-\delta_{\overline{r+s}} - \frac{\overline{r+s}}{T}-1}x_2^{\wt b - N - 1+i}(1+y)^{\wt u + N+ \delta_r + \frac{r}{T}}y^m\nn
&&\;\;\;\;\;\;\;\;\cdot Y_{\tilde{S}(W)}(Y((x_2+x_0)^{L(0)+ \delta_{\overline{r+s}} + \frac{\overline{r+s}}{T}}Y(u,y)v, x_0)b, x_2)w\nn
&=& \res_{x_0}\res_{x_2}\res_y \sum_{i =0}^Nx_0^{N-i-\delta_{\overline{r+s}} - \frac{\overline{r+s}}{T}-1}x_2^{\wt b - N - 1+i}(1+y)^{\wt u + N+ \delta_r + \frac{r}{T}}y^m\nn
&&\;\;\;\;\;\;\;\;\cdot (x_2+x_0)^{\wt u+\wt v+ \delta_{\overline{r+s}} + \frac{\overline{r+s}}{T}}Y_{\tilde{S}(W)}(Y(Y(u,(x_2+x_0)y)v, x_0)b, x_2)w\nn
&=& \res_{x_0}\res_{x_2}\res_{x_3} \sum_{i =0}^Nx_0^{N-i-\delta_{\overline{r+s}} - \frac{\overline{r+s}}{T}-1}x_2^{\wt b - N - 1+i}(x_2+x_0+x_3)^{\wt u + N+ \delta_r + \frac{r}{T}}x_3^m\nn
&&\;\;\;\;\;\;\;\;\cdot (x_2+x_0)^{\wt v-N+ \delta_{\overline{r+s}}+\frac{\overline{r+s}}{T}-\delta_r-\frac{r}{T} -m-1}Y_{\tilde{S}(W)}(Y(Y(u,x_3)v, x_0)b, x_2)w\nn
&=& \res_{x_0}\res_{x_2}\res_{x_1} \sum_{i =0}^Nx_0^{N-i-\delta_{\overline{r+s}} - \frac{\overline{r+s}}{T}-1}x_2^{\wt b - N - 1+i}(x_2+x_1)^{\wt u + N+ \delta_r + \frac{r}{T}}(x_1-x_0)^m\nn
&&\;\;\;\;\;\;\;\;\cdot (x_2+x_0)^{\wt v-N+ \delta_{\overline{r+s}}+\frac{\overline{r+s}}{T}-\delta_r-\frac{r}{T}-m-1}Y_{\tilde{S}(W)}(Y(u,x_1)Y(v, x_0)b, x_2)w\nn
&& -\; \res_{x_0}\res_{x_2}\res_{x_1} \sum_{i =0}^Nx_0^{N-i-\delta_{\overline{r+s}} - \frac{\overline{r+s}}{T}-1}x_2^{\wt b - N - 1+i}(x_2+x_1)^{\wt u + N+ \delta_r + \frac{r}{T}}(-x_0+x_1)^m\nn
&&\;\;\;\;\;\;\;\;\cdot (x_2+x_0)^{\wt v-N+ \delta_{\overline{r+s}}+\frac{\overline{r+s}}{T}-\delta_r-\frac{r}{T}-m-1}Y_{\tilde{S}(W)}(Y(v, x_0)Y(u,x_1)b, x_2)w\nn
\end{eqnarray*}

By checking the monomials involving $x_1$ in the first term of the
right-hand side of the equation above, the first term is an action of
elements in $O_g(V)$ on $w$, that is $0$ in $\tilde{S}(W)$.

For the second term, we shall check the monomials involving $x_0$. It is the action of an element of the form
\begin{equation}\label{newequation5}
\res_{x_0}(1+x_0)^{\wt v-N+ \delta_{\overline{r+s}}+\frac{\overline{r+s}}{T}-\delta_r-\frac{r}{T}-m-1}x_0^{m'}Y(v, x_0)\tilde{u}
\end{equation}
on $W$, where $m' \leq m + N-i-\delta_{\overline{r+s}} - \frac{\overline{r+s}}{T}-1$. Note that
\[
-N+ \delta_{\overline{r+s}}+\frac{\overline{r+s}}{T}-\delta_r-\frac{r}{T}-m-1 \geq \delta_s + \frac{s}{T} - 1,
\]
the element (\ref{newequation5}) can be written as an element of the form
\[
\res_{x_0}(1+x_0)^{\wt v+\delta_s + \frac{s}{T}-1}x_0^{m''}Y(v, x_0)\tilde{u},
\]
where $m'' = m' -N+ \delta_{\overline{r+s}}+\frac{\overline{r+s}}{T}-\delta_r-\frac{r}{T}-m-1 - \delta_s - \frac{s}{T} + 1 \leq -1-\delta_s$, hence the element (\ref{newequation5}) lies in $O_g(V)$ and the second term is $0$ in $\tilde{S}(W)$.  \epfv

Let $S_g(W) = \tilde{S}(W)/\mathcal{J}$. Then $S_g(W)$
is also a $\frac{1}{T}\N$-graded $T(V^g[t, t^{-1}])$-module. We can still
use elements of $T(V^g[t, t^{-1}])\otimes W$ to represent elements
of $S_g(W)$. But note that these elements now satisfy relations.
We equip $S_g(W)$ with the vertex operator map
\[
Y_{S_g(W)}: V \otimes S_g(W) \longrightarrow S_g(W)[[x^{\frac{1}{T}}, x^{-\frac{1}{T}}]]
\]
given by
\[
u \otimes w \rightarrow Y(u, x)w = Y_t(u, x)w.
\]

\begin{thm}
  The pair $(S_g(W), Y_{S_g(W)})$ is an admissible $g$-twisted
  $V$-module such that\\ $(S_g(W))_0 = W$.
\end{thm}
\pf As in $\tilde{S}(W)$, for $u \in V$ and $w \in S_g(W)$, we also
have $u(m)w = 0$ when $m > \wt u + \deg w - 1$. Clearly,
\[
Y_{S_g(W)}({\bf 1}, x) = I_{S_g(W)},
\]
where $I_{S_g(W)}$ is the identity operator on $S_g(W)$.

By Theorem \ref{main theorem}, where we specialize $k = \wt v + \deg w$ and $l = \wt u +
\deg w + \delta_r + \frac{r}{T}$ for $u \in V^{r}$, $S_g(W)$ satisfies weak associativity and hence is an admissible $g$-twisted $V$-module by Theorem \ref{main thm 1}. The claim that $(S_g(W))_0 = W$ follows from Lemma
\ref{lemma4}.\epfv

\begin{thm}
The functor $S_g$ has the following universal property:
Let $W$ be an $A_g(V)$-module. For any admissible twisted $V$-module
$\tilde{W}$ and any
$A_g(V)$-module map $f: W \to \Omega(\tilde{W})$, there exists a unique
$V$-module map $\tilde{f}: S_g(W)\to \tilde{W}$ such that
$\tilde{f}|_{W}=f$.
\end{thm}

In \cite{DLM2}, a functor, denoted by $\bar{M}$,
was constructed explicitly and was proved to satisfy the same universal property
above. The following result achieves our goal
of constructing this functor without dividing
relations corresponding to the commutator formula for weak modules:

\begin{cor}
The functor $S_g$ is equivalent to
the functor $\bar{M}$ constructed in \cite{DLM2}.
\end{cor}
\pf
This result follows immediately from the universal property.
\epfv

\def\refname{\hfil{REFERENCES}}

\noindent {\small \sc Department of Mathematics, University of Notre Dame,
278 Hurley Building, Notre Dame, IN 46556}
\vspace{1em}

\noindent {\em E-mail address}: jyang7@nd.edu
\end{document}